\documentclass[conference]{IEEEtran}
\IEEEoverridecommandlockouts
\usepackage{bm}
\usepackage{amsthm} 
\newtheorem{theorem}{Theorem}

\newtheorem{remark}{Remark}
\newtheorem{lemma}{Lemma}
\newtheorem{definition}{Definition}
\newtheorem{assumption}{Assumption}

\usepackage{algorithm}
\usepackage{algpseudocode} 
\usepackage{amsmath,amsfonts}
\allowdisplaybreaks[4]
\usepackage{array}
\usepackage[caption=false,font=normalsize,labelfont=sf,textfont=sf]{subfig}
\usepackage{textcomp}
\usepackage{stfloats}
\usepackage{url}
\usepackage{verbatim}
\usepackage{graphicx}
\usepackage{cite}
\usepackage{amssymb}
\usepackage{nomencl}
\makenomenclature
\usepackage{etoolbox}
\usepackage{url}
\usepackage{booktabs} 
\usepackage{xcolor}
\usepackage{multirow}
\newcommand{\xiny}{}
\newcommand{\xinyy}{}
\usepackage[prependcaption,colorinlistoftodos]{todonotes}

\newcommand{\icl}{}
\newcommand{\x}{}
\newcommand{\ic}{}

\def\BibTeX{{\rm B\kern-.05em{\sc i\kern-.025em b}\kern-.08em
    T\kern-.1667em\lower.7ex\hbox{E}\kern-.125emX}}

\makeatletter
\newcommand{\linebreakand}{%
  \end{@IEEEauthorhalign}
  \hfill\mbox{}\par
  \mbox{}\hfill\begin{@IEEEauthorhalign}
}
\makeatother

\begin{document}

\title{LiFT-MPC: Language-in-the-Loop Feedback Tuning of Cost Previews for MPC}

\author{Xinyi Yi and Ioannis Lestas%
\thanks{X. Yi and I. Lestas are with the Department of Engineering, University of Cambridge, Trumpington Street, Cambridge CB2 1PZ, UK.
Emails: {\tt\small <xy343, icl20>@cam.ac.uk}. 
Accepted for presentation at CDC 2026 and publication in the conference proceedings}
}

\addtolength{\topmargin}{0.3in}
\maketitle

\begin{abstract} 
In model predictive control (MPC) with time-varying objectives, \icl{predicted signals need to be often incorporated in the cost function, such as prices in energy system operation. These are, however, often difficult to predict from the historical trajectory of these signals alone, as they may depend on other contextual events. }
We propose \emph{LiFT-MPC}, \icl{an} MPC framework that integrates a \emph{LiFT} \icl{(Language-in-the-Loop Feedback Tuning)} correction \icl{scheme} 
to refine \icl{such predictions} 
within the MPC loop.
\icl{The prediction mechanism} is updated online via a \xiny{control-performance} loss \icl{function}, and we establish a performance guarantee for the resulting \icl{closed loop system.}
\icl{Numerical experiments using a realistic example of 
energy-storage management with real prices and}
news context \icl{to improve predictions,} demonstrate \icl{an improved  economic performance.} 
\end{abstract}

\section{Introduction}

In model predictive control (MPC) with time-varying objectives, preview quality strongly affects closed-loop performance. Inaccurate economic-signal \icl{predictions} 
may induce persistent input bias and degraded performance under receding-horizon control. Existing works have studied preview mismatch in MPC, including robustness and stability under disturbances and cost perturbations~\cite{berberich2025overview} and the effect of objective-term mismatch on suboptimality and practical stability~\cite{grune2020economic}. These results motivate \icl{the inclusion of preview refinements in MPC schemes}.

Preview signals are typically generated by numerical models, from classical time-series methods to recent learning-based approaches~\cite{lago2021forecasting,ansari2024chronos,rasul2024lagllama}. \x{However, even \ic{high quality} numerical previews may exhibit prediction errors when the previewed signal depends on contextual events that are difficult to infer from numerical history alone.} Electricity prices are a representative example, since outages and extreme weather can \x{strongly} affect price patterns, while information about such events is often conveyed through text. \x{In such settings, contextual information can provide complementary information~\cite{chattopadhyay2024context}, and \ic{the 
prediction error} can be treated as a residual to be corrected on top of a baseline numerical preview~\cite{junior2023hybrid}.}

\emph{Context-aware} MPC uses exogenous information to refine previews or adjust constraints. Many such methods remain \emph{context-then-control}: context is first used to improve prediction, and the resulting signal is then passed to MPC. Examples include event-triggered strategy switching~\cite{yang2023machine} and the use of semantic and geometric information in navigation to shape safety constraints and costs~\cite{10680398}. More recent \xiny{language-informed} methods move \ic{towards} a \emph{control-in-the-loop} paradigm, in which preview refinement is guided by downstream control performance~\cite{zeng2026llmmpc,wu2025instructmpc}.

Motivated by this viewpoint, we propose \emph{LiFT-MPC} (\textbf{L}anguage-\textbf{i}n-the-loop \textbf{F}eedback \textbf{T}uning MPC), \xiny{an} MPC framework that incorporates a \emph{LiFT} correction \icl{scheme} 
to refine \icl{the predictions of signals incorporated in the objective function of the MPC policy implemented.}
\x{A key feature of \emph{LiFT-MPC} is its residual correction structure: a baseline numerical preview is retained as the nominal MPC input, while contextual information is used only to generate an \ic{additive 
correction.}} 
\x{Unlike language-informed methods that directly generate disturbance previews from contextual inputs~\cite{zeng2026llmmpc,wu2025instructmpc}, \emph{LiFT-MPC} retains a baseline numerical preview and refines previewed economic signals that enter the MPC cost function. This makes the role of contextual information on the MPC explicit.}

The main contributions of this paper are threefold. First, we propose \emph{LiFT-MPC}, a language-\xiny{informed} preview-refinement framework for MPC with time-varying \icl{objectives.}
Second, we develop a delayed 
online update rule for \icl{the} \emph{LiFT} correction \icl{based on control-performance.} 
Third, we establish a performance guarantee for the resulting \icl{closed loop system that incorporates the online refinements.}
\xiny{Finally, we} validate \emph{LiFT-MPC} on \icl{a realistic example of} electricity-price-driven energy-storage \icl{management} \x{using real prices and news context in New South Wales, Australia, over 2015--2024, and evaluate its performance across numerical previews of different quality.}

The rest of the paper is organized as follows. Section~II presents the problem formulation. Section~III introduces the \emph{LiFT-MPC} framework and develops the associated theoretical results. Section~IV \ic{presents} the numerical experiments. Section~V concludes the paper. \icl{The proofs of the results have been included in the Appendix, which includes also more details on the numerical implementation of \ic{the} LiFT prediciton scheme.}

\section{Problem Formulation}\label{sec:preliminaries}
\subsection{Main notation}
The main notation is summarized in Table~\ref{tab:notation}. \xiny{We use $[T]:=\{0,1,\ldots,T-1\}$ to denote the full \xinyy{time} horizon \x{and $k$ to denote the preview-interval length}. We use $\{t,\ldots,\mathcal T(t)\}$ to denote the preview interval at time $t$, \xinyy{where $\mathcal T(t):=\min\{t+k-1,T-1\}$.} 
Throughout the paper, $\|\cdot\|$ denotes the Euclidean norm for vectors and the induced \(2\)-norm for matrices, unless otherwise specified. \icl{Sequences $x_t, x_{t+1}, \hdots, x_T$ where $t,T\in\mathbb{Z_+}$, $t\leq T$ will be denoted as $x_{t:T}$.} }

\begin{table}[htbp]
\centering
\caption{Main notation}
\label{tab:notation}
\small
\setlength{\tabcolsep}{4pt}
\renewcommand{\arraystretch}{1.08}
\begin{tabular}{p{1.55cm} p{6.35cm}}
\toprule
\textbf{Symbol} & \textbf{Description} \\
\midrule
$\xiny{T\in\mathbb{Z}_+}$ & Full-horizon length \\
$\xiny{k\in\mathbb{Z}_+}$ & Preview-interval length \\
$\xiny{\mathcal T(j)\in\mathbb{Z}_+}$ & \xiny{Terminal index of the preview interval starting at time $j$, with $\mathcal T(j):=\min\{j+k-1,T-1\}$} \\
\midrule
$\xiny{\pi_t\in\mathbb{R}^m}$ & \xiny{Real exogenous signal at time $t$} \\
$\xinyy{\pi_{i:\mathcal T(i)}}$ & \xiny{Real exogenous signal sequence \((\pi_i,\ldots,\pi_{\mathcal T(i)})\) over the realized interval $\{i,\ldots,\mathcal T(i)\}$} \\
\midrule
$\xiny{\hat{\pi}^{\mathrm{base}}_{\tau}\in\mathbb{R}^m}$ & \xiny{Baseline numerical preview at time $\tau$} \\
$\xinyy{\hat{\pi}^{\mathrm{base}}_{t:\mathcal T(t)}}$ & \xiny{Baseline numerical preview sequence \((\hat{\pi}^{\mathrm{base}}_t,\ldots,\hat{\pi}^{\mathrm{base}}_{\mathcal T(t)})\) over the preview interval $\{t,\ldots,\mathcal T(t)\}$} \\
$\xinyy{\Delta\pi_{t:\mathcal T(t)\mid t}}$ & \xiny{\emph{LiFT} correction sequence \((\Delta\pi_{t\mid t},\ldots,\Delta\pi_{\mathcal T(t)\mid t})\) \icl{generated at time $t$} over the preview interval.}\\
$\xiny{\hat{\pi}_{\tau\mid t}\in\mathbb{R}^m}$ & \xiny{\icl{Prediction generated at time $t$}} 
of the exogenous signal at time \icl{$\tau>t$}.\\  
$\xinyy{\hat{\pi}_{t:\mathcal T(t)\mid t}}$ & \xiny{
\icl{Sequence} \((\hat{\pi}_{t\mid t},\ldots,\hat{\pi}_{\mathcal T(t)\mid t})\).} \\ 
\midrule
$\xiny{s_{\tau\mid t}\in\mathbb{R}^n}$ & \xiny{Backward-recursion variable at time \icl{$\tau>t$ generated at time $t$} over the preview interval} \xinyy{(defined in (\ref{eq:s_recursion_mpc}))}\\
$\xinyy{s_{t:\mathcal T(t)\mid t}}$ & \xiny{
\icl{Sequence} \((s_{t\mid t},\ldots,s_{\mathcal T(t)\mid t})\)}\\ 
$\xiny{s_t^\star\in\mathbb{R}^n}$ & \xiny{Offline backward-recursion variable at time $t$, computed from the real exogenous signal sequence over the full horizon} \xinyy{(defined in (\ref{eq:s_recursion_mpcoff}))} \\
\midrule
$\xiny{s_{\tau\mid i}^{\mathrm{real}}\in\mathbb{R}^n}$ &\xiny{Backward-recursion variable \icl{at time $\qquad \qquad$ \hfill $\tau\in \{i,\ldots,\mathcal T(i)\}$}}  \xinyy{(defined in (\ref{realsignal}))} \\
$\xinyy{s_{i:\mathcal T(i)\mid i}^{\mathrm{real}}}$ & \xiny
\icl{Sequence \((s_{i\mid i}^{\mathrm{real}},\ldots,s_{\mathcal T(i)\mid i}^{\mathrm{real}})\)} 
\\
\bottomrule
\end{tabular}
\end{table}

\subsection{Problem Setup}
\subsubsection{System dynamics and time-varying objective}
Consider the linear system 
\begin{equation}
x_{t+1}=A x_t + B u_t,\quad t\in[T]:=\{0,1,\ldots,T-1\},
\end{equation}
\xinyy{where $x_t \in \mathbb{R}^n$ and $u_t \in \mathbb{R}^m$.}
The finite-horizon cost is defined as:
\begin{subequations}\label{eq:fh_cost_general}
\begin{align}
J_T
:=\;&
\sum_{t=0}^{T-1}\ell_t(x_t,u_t;\pi_t) + V_f(x_T),\label{eq:fh_cost_generala}\\
\ell_t(x_t,u_t;\pi_t)
:=\;&
x_t^\top Q x_t + u_t^\top R u_t + 2\alpha_\pi\,\pi_t^\top u_t,\label{eq:stage_cost}
\end{align}
\end{subequations}
where $Q\succeq 0$, $R\succ 0$, and $\alpha_\pi>0$.

\begin{assumption}\label{ass:system}
$(A,B)$ is controllable and $(Q^{1/2},A)$ is detectable.
\end{assumption}

Under Assumption~\ref{ass:system}, let $P\succeq0$ denote the stabilizing solution of the discrete-time algebraic Riccati equation
\begin{equation}\label{eq:dare}
P
=
Q+A^\top P A
-
A^\top P B\left(R+B^\top P B\right)^{-1}B^\top P A .
\end{equation}
We choose the terminal cost\footnote{This terminal choice provides a standard approximation of the tail cost beyond the finite optimization interval, helps prevent myopic end-of-horizon behavior, and yields the affine MPC law used later.} as $V_f(x_T)=x_T^\top P x_T$.

\begin{assumption}\label{ass:price_bounds}
The exogenous signal \xinyy{$\pi_t$} satisfies $\underline{\pi} \le \pi_t \le \overline{\pi}$ for all $t$, where $\underline{\pi},\overline{\pi}\in\mathbb{R}^m$ are given bounds \xiny{for $\pi_t$}.
\end{assumption}

\subsubsection{Offline benchmark (Off-OPT)}
Suppose the controller has access to the real sequence $\pi_{0:T-1}$. Define
\begin{equation}\label{eq:pre_opt}
\begin{aligned}
u^{\mathrm{Off\text{-}OPT}}_{0:T-1}
:=
&\arg\min_{\icl{u_{0:T-1}}}\;
\sum_{t=0}^{T-1}\ell_t(x_t,u_t;\pi_t) + V_f(x_T)\\
\text{s.t.}\;
&x_{t+1}=A x_t + B u_t,\qquad t=0,\ldots,T-1 .
\end{aligned}
\end{equation}
Using dynamic programming, let
$V_T(x_T)=V_f(x_T)=x_T^\top P x_T$ and
$V_t(x_t)=\min_{u_t}\bigl(\ell_t(x_t,u_t;\pi_t)+V_{t+1}(A x_t+B u_t)\bigr)$, $t=0,\ldots,T-1$.
Consider the affine-quadratic ansatz
\xiny{$V_t(x_t)=x_t^\top P x_t+2s_t^\top x_t$},
where $P$ is the stabilizing solution of \eqref{eq:dare}, $s_T=0$.
Substituting this ansatz into the Bellman recursion gives
\xiny{\begin{equation}\label{eq:recursionansatz}
\begin{aligned}
V_t(x_t)=\min_{u_t}\Big[
&x_t^\top Qx_t+u_t^\top Ru_t+2\alpha_\pi \pi_t^\top u_t\\
&+(Ax_t+Bu_t)^\top P(Ax_t+Bu_t)\\
&+2s_{t+1}^\top(Ax_t+Bu_t)
\Big].
\end{aligned}
\end{equation}}
Define
$H:=R+B^\top P B\succ0$ and $K:=H^{-1}B^\top P A$.
Collecting the $u_t$-dependent terms in \eqref{eq:recursionansatz} and applying the first-order optimality condition yield the optimal affine control law
\begin{equation}\label{offoptlaw}
u_t^{\mathrm{Off\text{-}OPT}}
=
-Kx_t-H^{-1}\!\bigl(B^\top s_{t+1}^\star+\alpha_\pi\pi_t\bigr).
\end{equation}
Substituting \eqref{offoptlaw} back into \eqref{eq:recursionansatz} and comparing the coefficients of the linear term in $x_t$ gives the backward recursion
\begin{equation}\label{eq:s_recursion_mpcoff}
s_t^\star
=
(A-BK)^\top s_{t+1}^\star-\alpha_\pi K^\top \pi_t,
\qquad s_T^\star=0.
\end{equation}

\begin{definition}\label{defphi}
Let $\Phi:=A-BK$. Since $P$ is the stabilizing solution of \eqref{eq:dare}, $\Phi$ is asymptotically stable, \xiny{(i.e. its spectral radius is less than 1)}. Hence there exist constants $c_\Phi>0$ and $\rho\in(0,1)$ such that $\|(\Phi^\top)^j\|\le c_\Phi\rho^j$ for all $j\ge0$. \xiny{Here $(\Phi^\top)^j$ denotes the $j$-th matrix power of $\Phi^\top$, with $(\Phi^\top)^0:=I$.}
\end{definition}

\xiny{Recursion \eqref{eq:s_recursion_mpcoff} yields}
\begin{equation}\label{recurcal}
s_{t+1}^\star
=
-\sum_{j=0}^{T-t-2}(\Phi^\top)^j \alpha_\pi K^\top \pi_{t+1+j}, t=0,\ldots,T-2.
\end{equation}
The associated closed-loop performance is
\begin{equation}\label{eq:J_preopt_def}
J^{\mathrm{Off\text{-}OPT}}
:=
\sum_{t=0}^{T-1}\ell_t\!\bigl(x_t,u_t^{\mathrm{Off\text{-}OPT}};\pi_t\bigr)+V_f(x_T).
\end{equation}

\subsection{Preview-based MPC}\label{subsec:mpc}
Because the stage cost \eqref{eq:stage_cost} depends on the time-varying signal $\pi_t$, whose future values 
\icl{can only be predicted,} 
we adopt a receding-horizon implementation. Fix a preview interval length $k\in\mathbb N$ and define $\mathcal T(t):=\min\{t+k-1,T-1\}$. At time $t\in[T]$, let
$\hat{\pi}_{t:\mathcal T(t)\mid t}:=(\hat{\pi}_{t\mid t},\hat{\pi}_{t+1\mid t},\ldots,\hat{\pi}_{\mathcal T(t)\mid t})$
denote the exogenous signal preview \icl{(prediction)} formed at time $t$ over the preview interval $t,\ldots,\mathcal T(t)$.

Given the current state $x_t$, MPC solves
\begin{equation}\label{eq:mpc_problem}
\begin{aligned}
u^{\mathrm{MPC}}_{t:\mathcal T(t)}
:=&
\arg\min_{u_{t:\mathcal T(t)}}\;
\sum_{\tau=t}^{\mathcal T(t)}
\ell_\tau\!\bigl(x_\tau,u_\tau;\hat{\pi}_{\tau\mid t}\bigr)
+V_f\!\bigl(x_{\mathcal T(t)+1}\bigr)\\
&\text{s.t.}\;
x_{\tau+1}=A x_\tau + B u_\tau, t\le \tau\le \mathcal T(t).
\end{aligned}
\end{equation}
Let $u^\star_{\tau\mid t}$ and $x_{\tau\mid t}$ denote the optimal input and predicted state at time $\tau$ along the solution of \eqref{eq:mpc_problem}. By the same dynamic-programming argument as in the offline benchmark, \icl{we consider the backward recursion}
\begin{equation}\label{eq:s_recursion_mpc}
s_{\tau\mid t}
=
(A-BK)^\top s_{\tau+1\mid t}
-\alpha_\pi K^\top \hat{\pi}_{\tau\mid t},
s_{\mathcal T(t)+1\mid t}=0.
\end{equation}
The corresponding affine control law is 
\begin{equation}\notag
u_{\tau\mid t}^\star
=
-Kx_{\tau\mid t}
-H^{-1}\!\Bigl(B^\top s_{\tau+1\mid t}+\alpha_\pi\hat{\pi}_{\tau\mid t}\Bigr).
\end{equation}
Receding-horizon MPC applies the first control input,
\begin{equation}\label{actmpc}
u_t^{\mathrm{MPC}}
:=
u_{t\mid t}^\star
=
-Kx_t-H^{-1}\!\Bigl(B^\top s_{t+1\mid t}+\alpha_\pi\hat{\pi}_{t\mid t}\Bigr),
\end{equation}
and \icl{recalculates the input} 
at time $t+1$ with a new preview.

For $t=0,\ldots,\mathcal T(t)-1$, \xiny{\eqref{eq:s_recursion_mpc} gives}
\begin{equation}\label{mpcs}
s_{t+1\mid t}
=
-\sum_{j=0}^{\mathcal T(t)-t-1}(\Phi^\top)^j \alpha_\pi K^\top \hat{\pi}_{t+1+j\mid t}.
\end{equation}
The resulting closed-loop performance \icl{evaluated} under the real \xiny{signal $\pi_t$ \icl{is given by the expression}}:
\begin{equation}\label{eq:J_mpc_def}
J^{\mathrm{MPC}}
:=
\sum_{t=0}^{T-1}\ell_t\!\bigl(x_t,u_t^{\mathrm{MPC}};\pi_t\bigr)+V_f(x_T).
\end{equation}
\icl{Note that} in \eqref{eq:J_mpc_def} \icl{the control input is determined} 
using the preview $\hat{\pi}_{t:\mathcal T(t)\mid t}$, while \icl{the performance metric specified} is evaluated using the real signal $\pi_t$.

\begin{remark}[Regularity of the stage cost]\label{rem:cost_regularity}
For all $t$ and all $\pi,\pi'\in\{\xi\in\mathbb{R}^m:\underline{\pi}\le \xi\le \overline{\pi}\}$, the stage cost $\ell_t(x,u;\pi)$ is continuously differentiable and strictly convex in $u$. Moreover, for any $\bar u>0$ and any $u$ satisfying $\|u\|\le \bar u$,
\[
\bigl|\ell_t(x,u;\pi)-\ell_t(x,u;\pi')\bigr|
\le
2\alpha_\pi\,\bar u\,\|\pi-\pi'\|.
\]
\end{remark}

\subsection{Control-relevant optimality gap}\label{subsec:gap}

\subsubsection{Exact optimality gap}
Here $s_{t+1\mid t}$ denotes the preview backward-recursion variable generated by \xiny{(\ref{eq:s_recursion_mpc}), (\ref{mpcs})} from the preview $\hat{\pi}_{t:\mathcal T(t)\mid t}$, while $s_t^\star$ denotes the backward-recursion variable generated by \xiny{(\ref{eq:s_recursion_mpcoff}), (\ref{recurcal})} from the real signal sequence over the full horizon.

\begin{lemma}\label{lem:quadratic_gap}
The optimality gap \xinyy{satisfies the relation}
\begin{equation}\label{eq:gap_identity_suma}
J^{\mathrm{MPC}}-J^{\mathrm{Off\text{-}OPT}}
=
\sum_{t=0}^{T-1}\psi_t^\top H\psi_t
=
\sum_{t=0}^{T-1}F_t,
\end{equation}
where \xiny{$\psi_t$} and the per-round loss $F_t$ are defined as
\begin{subequations}\label{eq:gap_identity_sumover}
\begin{align}
\psi_t
&:=
-H^{-1}\Bigl(
B^\top\bigl(s_{t+1\mid t}-s_{t+1}^\star\bigr)
+\alpha_\pi\bigl(\hat{\pi}_{t\mid t}-\pi_t\bigr)
\Bigr),\label{eq:gap_identity_sumc}\\
F_t&:=\psi_t^\top H\psi_t.\label{eq:gap_identity_sumb}
\end{align}
\end{subequations}
\end{lemma}
Lemma~\ref{lem:quadratic_gap} quantifies the loss induced by preview mismatch, \icl{i.e. the performance loss due to the determination of the control policy using the predicted $\hat \pi_t$ rather than its actual value}. \icl{Note, however, that} the exact per-round loss $F_t$ depends on $s_{t+1}^\star$, which is generated from the real signal over the full horizon and is therefore unavailable online.

\subsubsection{Online computable optimality gap}
At time $t\ge k$, let $i,\ldots,\mathcal T(i)$ denote the latest realized interval, \icl{i.e.} $\mathcal T(i)+1=t$. Using the revealed real signals over this realized interval, 
\begin{equation}\label{realsignal}
\begin{aligned}
&s^{\mathrm{real}}_{\tau\mid i}
=
(A-BK)^\top s^{\mathrm{real}}_{\tau+1\mid i}
-\alpha_\pi K^\top \pi_\tau, \tau=i,\ldots,\mathcal T(i),\\
&s^{\mathrm{real}}_{\mathcal T(i)+1\mid i}=0.
\end{aligned}
\end{equation}
\xiny{Expanding} \eqref{realsignal} gives
\begin{equation}\label{reals}
s_{i+1\mid i}^{\mathrm{real}}
=
-\sum_{j=0}^{\mathcal T(i)-i-1}(\Phi^\top)^j\alpha_\pi K^\top \pi_{i+1+j}.
\end{equation}

\begin{definition}\label{prop:surrogate_loss}
Define the 
\icl{proxy per-round loss}
$L_i:=\widehat\psi_i^\top H\widehat\psi_i$, \x{where $\widehat\psi_i
:=
-H^{-1}\Bigl(
B^\top\bigl(s_{i+1\mid i}-s_{i+1\mid i}^{\mathrm{real}}\bigr)
+\alpha_\pi\bigl(\hat{\pi}_{i\mid i}-\pi_i\bigr)
\Bigr)$.}
\end{definition}

Comparing \eqref{reals} with \eqref{recurcal} yields
\begin{equation}\notag
s_{i+1}^\star-s_{i+1\mid i}^{\mathrm{real}}
=
-\sum_{j=\mathcal T(i)-i}^{T-i-2}(\Phi^\top)^j \alpha_\pi K^\top \pi_{i+1+j},
\end{equation}
and therefore
\begin{equation}\label{eq:s_trunc_bound_final_polished_appendix}
\bigl\|s_{i+1\mid i}^{\mathrm{real}}-s_{i+1}^\star\bigr\|
\le
\alpha_\pi \|K\|\,\bar\pi\,
\frac{c_\Phi\,\rho^{\mathcal T(i)-i}}{1-\rho}.
\end{equation} 

\begin{remark}
   \icl{ It should be noted that \x{$L_i$} defined above quantifies delayed performance loss, it is, however, information available at time $t$ to the controller. In the next section this will be used as a basis for refining the prediction mechanism for $\pi_t$.} 
\end{remark}

\section{\emph{LiFT-MPC} Framework}
Baseline numerical previews capture the dominant temporal structure of the exogenous signal. \x{However, prior work suggests that even strong baseline previews can be improved by incorporating additional contextual information ~\cite{chattopadhyay2024context} (e.g., relevant background or auxiliary information beyond the baseline numerical model’s standard inputs) and by residual correction (i.e., correcting the remaining error not captured by the baseline numerical model)~\cite{junior2023hybrid}. Motivated by this idea, we retain the numerical preview as the nominal input to MPC and use the \emph{LiFT} module only to generate an additive correction. Unlike methods that modify the predictor itself~\cite{chattopadhyay2024context} or learn a separate residual model from numerical history~\cite{junior2023hybrid}, our approach keeps the baseline numerical preview $\hat{\pi}^{\mathrm{base}}_{t:\mathcal T(t)}$ fixed and applies a 
context-dependent residual correction $\Delta\pi_{t:\mathcal T(t)\mid t}$.}

\ic{In particular, the corrected preview $\hat{\pi}_{t:\mathcal T(t)\mid t}$ used by MPC is defined as the sum of a numerical preview and residual correction:}
\begin{equation}\label{eq:surrogate_preview}
\hat{\pi}_{t:\mathcal T(t)\mid t}
=
\hat{\pi}^{\mathrm{base}}_{t:\mathcal T(t)}
+
\Delta\pi_{t:\mathcal T(t)\mid t}.
\end{equation}
\x{This makes the effect of contextual information on the control policy explicit. Throughout the subsequent analysis, the baseline numerical preview $\hat{\pi}^{\mathrm{base}}_{t:\mathcal T(t)}$ is treated as fixed, and all comparisons are made relative to this same baseline preview.}

\subsection{Residual correction}\label{liftcons}
\x{To incorporate contextual information in a structured way, we represent the residual correction as a weighted combination of a finite collection of fixed historical correction trajectories associated with different operating scenarios (e.g., equipment outages and extreme weather conditions in power systems) \ic{This is quantified in equation \eqref{mixweight} later in this section}.}

\begin{assumption}[Fixed context over the full horizon]\label{ass:fixed_context_allocation}
\x{The contextual information used by the \emph{LiFT} module is collected before the full horizon $[T]$ starts\footnote{This choice is consistent with energy-system operation settings in which the full horizon corresponds to one operating day and the contextual information is observed before the control process over that horizon begins, remaining unchanged throughout it~\cite{antoniadou2022scenario,liu2015bidding}.} and is treated as fixed over the full \ic{horizon.}} 
\end{assumption}

\subsubsection{\emph{LiFT} scenario-mixture correction}
\x{Let $\mathcal S$ denote a finite scenario set. For each $s\in\mathcal S$, let $h^s_{0:T-1}$ denote a fixed historical correction trajectory over the full horizon, constructed from historical data\footnote{The construction of the historical correction trajectories is deferred to Appendix~E. Other construction procedures are also possible.}, and let $h^s_{t:\mathcal T(t)} := (h^s_t,\ldots,h^s_{\mathcal T(t)})$ denote the truncated trajectory over the preview \ic{interval.}} 

\x{Given the fixed contextual information, an initial scenario-weight vector $\boldsymbol w$ is computed. \ic{The specific way this is constructed} is not needed for the theoretical development and is deferred to Appendix~F and the numerical experiments. 
To better exploit information revealed along the closed-loop trajectory, we introduce a refinement parameter $\theta\in\Theta$ to \ic{refine $\boldsymbol w$ .} 
The refined scenario-weight vector is then defined \ic{as}}
\begin{equation}\label{adaptweight}
\boldsymbol w_{\theta}=f(\theta,\boldsymbol w).
\end{equation}
\x{\ic{where $\boldsymbol w_{\theta}:=[w_\theta^s]_{s\in\mathcal S}$, and $w_\theta^s\in\mathbb R$ denotes} the weight assigned to scenario $s$ under the fixed contextual \ic{information $\boldsymbol w$} and the refinement parameter $\theta$. The residual correction in \eqref{eq:surrogate_preview} over the preview interval is then defined as}
\begin{equation}\label{mixweight}
\Delta\pi_{t:\mathcal T(t)\mid t}(\theta)
=
\sum_{s\in\mathcal S} w_\theta^s\, h^s_{t:\mathcal T(t)}.
\end{equation}
\x{At time $t$, the controller uses the current refinement parameter $\theta_t$ and forms the corrected preview (\ref{eq:surrogate_preview}) as $\hat{\pi}_{t:\mathcal T(t)\mid t}
=
\hat{\pi}^{\mathrm{base}}_{t:\mathcal T(t)}
+
\Delta\pi_{t:\mathcal T(t)\mid t}(\theta_t)$.} 

\begin{assumption}\label{assumeaffine}
\x{The map $\theta\mapsto f\bigl(\theta,\boldsymbol w\bigr)$ in (\ref{adaptweight}) is differentiable and affine in $\theta$.}
\end{assumption}

\subsubsection{Delayed update rule for the \emph{LiFT} parameter}
\x{At time $t\ge k$, let $i,\ldots,\mathcal T(i)$ denote the latest realized interval, \icl{i.e.} $\mathcal T(i)+1=t$. Over the realized interval, $\hat{\pi}_{i:\mathcal T(i)\mid i}(\theta_i)$ and $s_{i+1\mid i}$ used by the MPC can be calculated as:} $\hat{\pi}_{i:\mathcal T(i)\mid i}(\theta_i)
=
\hat{\pi}^{\mathrm{base}}_{i:\mathcal T(i)}
+
\Delta\pi_{i:\mathcal T(i)\mid i}(\theta_i)$,
and by \eqref{mpcs}, $s_{i+1\mid i}(\theta_i)
=
-\sum_{j=0}^{\mathcal T(i)-i-1}(\Phi^\top)^j\alpha_\pi K^\top
\hat{\pi}_{i+1+j\mid i}(\theta_i)$.
\x{The proxy per-round loss $L_i(\theta_i)$ in Definition \ref{prop:surrogate_loss} can be calculated as:}\footnote{Other loss functions can also be used, provided they act as tractable proxies for the exact gap $\psi_t$ and the exact per-round loss $F_t$ in \eqref{eq:gap_identity_sumover}. In the performance bound derived later, this choice affects only the surrogate-mismatch term in \eqref{eq:regret_bound_general}.}
\begin{equation}\label{thetacompu}
L_i(\theta_i):=
\widehat\psi_i(\theta_i)^\top H\,\widehat\psi_i(\theta_i),
\end{equation}
where $\widehat\psi_i(\theta_i)
:=-H^{-1}\Bigl(
B^\top\bigl(s_{i+1\mid i}(\theta_i)-s_{i+1\mid i}^{\mathrm{real}}\bigr)
+\alpha_\pi\bigl(\hat{\pi}_{i\mid i}(\theta_i)-\pi_i\bigr)
\Bigr)$.

\begin{assumption}[\x{LiFT parameter set and bounded gradients}]\label{ass:lift_regularity}
Let $\Theta\subset\mathbb R^d$ be convex and compact, and suppose $\theta_t\in\Theta$ for all $t\in[T]$, with \x{$\sup_{\theta,\theta'\in\Theta}\|\theta-\theta'\|\le D$}. Assume also that $\|\nabla_\theta L_t(\theta)\|\le G$ for all \x{$t\in[T]$} and $\theta\in\Theta$.
\end{assumption}

\xiny{Using the proxy per-round loss in \eqref{thetacompu}}, we apply delayed gradient descent with a non-increasing stepsize $\eta_t$:
\begin{equation}\label{eq:lift_update}
\theta_{t+1}
=
\theta_t-\eta_t\nabla_\theta L_i(\theta_i),t\ge k,
\end{equation}
\xiny{where $\nabla_\theta L_i(\theta_i)$ denotes} the gradient of the loss associated with time $i$, evaluated at the parameter value $\theta_i$ used at that time. For $t\le k-1$, we keep $\theta_t=\theta_0$.

\subsection{Performance bound with delayed updates}\label{subsec:regret}

\subsubsection{Loss discrepancy between exact and proxy per-round losses}
Analogously to \eqref{thetacompu}, over the
realized interval $i,\ldots,\mathcal T(i)$, \x{the exact per-round loss in Lemma \ref{lem:quadratic_gap} is calculated as:}
\begin{equation}\label{thetaexat}
F_i(\theta_i)
:=
\psi_i(\theta_i)^\top H\psi_i(\theta_i),
\end{equation}
where $\psi_i(\theta_i)
:=
-H^{-1}\Bigl(
B^\top\bigl(s_{i+1\mid i}(\theta_i)-s_{i+1}^\star\bigr)+\alpha_\pi\bigl(\hat{\pi}_{i\mid i}(\theta_i)-\pi_i\bigr)
\Bigr)$.
\begin{definition}\label{rhole}
Let $R_i(\theta):=\nabla_\theta\psi_i(\theta)$. Under Assumption~\ref{assumeaffine}, $\psi_i(\theta)$ is affine in $\theta$, so $R_i(\theta)$ is independent of $\theta$ for each fixed $i$. Since the horizon is finite, there exists a constant $\bar R>0$ such that $\sup_{\theta\in\Theta}\|R_i(\theta)\|\le \bar R$ for all $i\in[T]$.
\end{definition}

\begin{definition}\label{def:loss_discrepancy}
For differentiable losses $L_1,L_2:\Theta\to\mathbb R$, define
$\mathrm{LD}(L_1,L_2)
:=
\sup_{\theta\in\Theta}
\bigl\|\nabla_\theta L_1(\theta)-\nabla_\theta L_2(\theta)\bigr\|$.
\end{definition}

\begin{lemma}\label{lem:LD_bound_final}
For every $i\in[T]$, the proxy per-round loss $L_i(\theta_i)$ in \eqref{thetacompu} and the exact per-round loss $F_i(\theta_i)$ in \eqref{thetaexat} satisfy
\begin{equation}\label{eq:LD_bound_explicit_polished}
\mathrm{LD}(L_i,F_i)
\le
C_{\mathrm{LD}}\rho^{\mathcal T(i)-i},
\end{equation}
where $C_{\mathrm{LD}}
:=
2\bar R\|H\|\,\|B\|\,\|H^{-1}\|\,
\alpha_\pi\|K\|\,\bar\pi\,\frac{c_\Phi}{1-\rho}$,
\xiny{and $c_\Phi$ and $\rho$ are defined in Definition~\ref{defphi}.}
In particular, when $\mathcal T(i)=i+k-1$, the above bound becomes \(O(\rho^{k-1})\).
\end{lemma}

\subsubsection{Best fixed refinement parameter (\emph{LiFT}-OPT)}
\x{Under Assumption~\ref{ass:fixed_context_allocation}, the contextual information is fixed over the full horizon, and the online refinement updates the scenario weights $\boldsymbol w_{\theta}$ in \eqref{adaptweight} as estimates of a fixed underlying scenario weighting over the horizon. Accordingly, 
\ic{a natural} benchmark is the best fixed refinement parameter over the full \ic{horizon.}}

\x{Let $\theta_{0:T-1}:=(\theta_0,\ldots,\theta_{T-1})$ denote the online parameter sequence generated by \eqref{eq:lift_update}. The realized closed-loop cost (\ref{eq:J_mpc_def}) is defined as:}
\begin{equation}
J^{\mathrm{MPC}}(\theta_{0:T-1})
:=
\sum_{t=0}^{T-1}\ell_t\!\bigl(x_t,u_t^{\mathrm{MPC}}(\theta_t);\pi_t\bigr)+V_f(x_T).
\end{equation}
\x{For any fixed $\theta\in\Theta$, let $J^{\mathrm{MPC}}(\theta)
:=
\sum_{t=0}^{T-1}\ell_t\!\bigl(x_t,u_t^{\mathrm{MPC}}(\theta);\pi_t\bigr)+V_f(x_T)$ denote the corresponding realized closed-loop cost when the same refinement parameter $\theta$ is used over the full horizon. Define the best fixed refinement parameter $\theta^\star$ by}
\begin{equation}
\theta^\star\in\arg\min_{\theta\in\Theta}J^{\mathrm{MPC}}(\theta),
\qquad
J^{\mathrm{\emph{LiFT}\text{-}OPT}}:=J^{\mathrm{MPC}}(\theta^\star).
\end{equation}
\x{Thus, $J^{\mathrm{\emph{LiFT}\text{-}OPT}}$ is the best realized closed-loop cost achievable within the proposed \emph{LiFT} refinement structure when the same fixed refinement parameter is used over the full horizon.}

\begin{theorem}[Performance bound with delayed updates]
\label{thm:regret_delay_ld}
Fix a preview interval length $k\in\mathbb N$ and consider the delayed update \eqref{eq:lift_update}. Under Assumptions~\ref{ass:system}--\ref{ass:lift_regularity},
\begin{equation}\label{eq:regret_bound_general}
\begin{aligned}
&J^{\mathrm{MPC}}(\theta_{0:T-1})-J^{\mathrm{\emph{LiFT}\text{-}OPT}}\\
\le\;&
\underbrace{\frac{D^2}{2\eta_{T-1}}
+\frac{G^2}{2}\sum_{t=0}^{T-1}\eta_t}_{\text{standard online-gradient term}}+\underbrace{(k-1)\!\left(GD+G^2\sum_{t=0}^{T-1}\eta_t\right)}_{\text{delay term}}\\
&+\underbrace{D\sum_{t=0}^{T-1}C_{\mathrm{LD}}\rho^{\mathcal T(t)-t}}_{\text{surrogate-mismatch term}}.
\end{aligned}
\end{equation}
Choosing $\eta_t=\frac{D}{G\sqrt{t+1}}$ yields
$\frac{D^2}{2\eta_{T-1}} +\frac{G^2}{2}\sum_{t=0}^{T-1}\eta_t=2GD\sqrt{T}$.
\end{theorem}
\x{Theorem~\ref{thm:regret_delay_ld} bounds the realized closed-loop cost of the online refinement sequence relative to \emph{LiFT}-OPT, which is the best realizable closed-loop cost within the proposed \emph{LiFT} refinement structure over the full horizon. The bound separates the effects of online refinement, delayed feedback, and surrogate mismatch of loss functions in the MPC loop.}

\subsection{Practical tracking stability under preview mismatch}\label{subsec:stability_short}

We now state a practical tracking guarantee for the closed-loop system induced by preview-based MPC. Under the affine MPC law \eqref{actmpc}, the implemented input can be written as
\begin{equation}\label{eq:u_affine_split}
u_t=-Kx_t+v_t,
v_t:=-H^{-1}\!\Bigl(B^\top s_{t+1\mid t}(\theta_t)+\alpha_\pi\hat{\pi}_{t\mid t}(\theta_t)\Bigr).
\end{equation}
Let $(x_t^\star,u_t^\star)$ denote the trajectory-input pair induced by the offline benchmark \eqref{eq:pre_opt}. Equivalently, under the affine law \eqref{offoptlaw},
\begin{equation}\label{eq:u_star_affine_split}
u_t^\star=-Kx_t^\star+v_t^\star,
v_t^\star:=-H^{-1}\!\Bigl(B^\top s_{t+1}^\star+\alpha_\pi\pi_t\Bigr).
\end{equation}
Define the feedforward mismatch by $\delta_t:=v_t-v_t^\star$, and assume that \xiny{$\|\delta_t\|\le\bar\delta$ for all $t$ \xinyy{in the time horizon $[T]$}, for some constant $\bar\delta>0$.}

\begin{theorem}[Practical tracking stability]\label{thm:practical_stab_short}
Under the notation of Definition~\ref{defphi},
$\|x_t-x_t^\star\|
\le
c_\Phi\rho^t\|x_0-x_0^\star\|
+
\frac{c_\Phi\|B\|}{1-\rho}\bar\delta,
\forall t\ge0$.
Hence the tracking error converges exponentially to a tube whose radius scales linearly with $\bar\delta$.
\end{theorem}

\section{NUMERICAL EXPERIMENT}
\x{We evaluate whether day-ahead news context can improve the closed-loop performance $J^{\mathrm{MPC}}$ in \eqref{eq:J_mpc_def} for an MPC-operated battery energy storage system (ESS) using the proposed \emph{LiFT} correction module and real electricity-price and news data from New South Wales, Australia, over 2015--2024~\cite{bi2025nswepnews}. The previewed signal is the 48-step intra-day electricity-price trajectory.}

\subsection{Energy storage system (ESS) operation setting}
\label{subsec:ess_problem_formulation}

We consider an intra-day ESS operation problem over one day with 48 half-hour intervals.

\subsubsection{ESS operation model}
Let $x_t^0\in[0,E_{\max}]$ denote the state of charge (SoC), with $E_{\max}=200$ MWh, and let $u_t\in[-P_{\max},P_{\max}]$ denote the charging/discharging power, with $P_{\max}=100$ MW\footnote{This scale is comparable to utility-scale battery installations such as the Terang BESS in Australia and the Derrymeen BESS in Northern Ireland.}, where $u_t>0$ and $u_t<0$ correspond to charging and discharging, respectively.\footnote{In the numerical implementation, the executed control input is clipped to satisfy the power and SoC limits. \x{Since the theory in Section~III is developed for the unclipped affine policy, the effect of clipping on the performance bounds is not analyzed here and is left for future work.}} 
We adopt the energy-balance model $x_{t+1}^0=x_t^0+\Delta u_t$, \x{with} $\Delta=\tfrac12$ hours. To match the formulation in Section~II, we use the shifted state $x_t:=x_t^0-x_{\mathrm{ref}}$, where $x_{\mathrm{ref}}=\tfrac12 E_{\max}$, so that \x{$x_{t+1}=x_t+\Delta u_t$.}

The stage cost is $\ell_t(x_t,u_t;\pi_t)=Qx_t^2+Ru_t^2+2\alpha_\pi \pi_t u_t \Delta$, where $\pi_t\in\mathbb R$ denotes the electricity price \x{at time $t$}, \x{$Q=10^{-4}$ and $R=10^{-3}$ provide mild SoC regularization and input smoothing, and $\alpha_\pi=2$ captures the arbitrage incentive.} The terminal cost is \x{$V_f(x_T)=x_T^\top P x_T$}, where $P$ is the stabilizing solution of \eqref{eq:dare}\x{.} The economic reward over \x{one} horizon is defined as $\textstyle
\sum_{t=0}^{47} -\pi_t u_t \Delta$,
which corresponds to arbitrage profit under real prices. \x{The primary closed-loop performance is evaluated by $J^{\mathrm{MPC}}$ in \eqref{eq:J_mpc_def}.}

\subsubsection{Numerical baseline preview}
In our setting, a day-ahead numerical price preview provides the baseline signal used by the controller, while day-ahead news text provides contextual information for refining that preview. \x{This is standard in energy-system operation~\cite{antoniadou2022scenario,liu2015bidding}.} \x{Baseline preview models are trained on 2015--2017 and evaluated on 2018--2022 for model selection. Among the candidates, AR(48)~\cite{lago2021forecasting} and Chronos-T5~\cite{ansari2024chronos} are used as representative baselines, achieving the worst and best preview accuracy, respectively, on the 2018--2022 evaluation period.}

\subsubsection{Scenario library}
For each baseline \x{numerical} preview, we construct a separate \emph{LiFT} scenario library from prices and news over 2018--2022, so that the \x{fixed historical correction trajectories $h^s_{0:T-1}$} \xiny{are matched to the corresponding baseline}. \xiny{Details are deferred to Appendix~E.} 

\subsubsection{Preview-based MPC}
We apply the preview-based MPC formulation in Section~II with preview length $k=8$ (4 hours). At each time $t$, the controller uses the preview $\hat{\pi}_{t:\mathcal T(t)\mid t}$ and implements the input \eqref{actmpc}. The preview is constructed via \eqref{eq:surrogate_preview}--\eqref{mixweight}, where the \emph{LiFT} module refines the baseline numerical \xiny{preview} using day-ahead news.

\subsubsection{\emph{LiFT} correction in experiments}
\label{subsec:price_prediction}
\x{In the numerical experiments, the contextual information is given by day-ahead news and is fixed over the full intra-day horizon, consistent with Section~\ref{liftcons}. The initial scenario-weight vector $\boldsymbol w$ is constructed from this news using the pretrained encoder \texttt{sentence-transformers/all-MiniLM-L6-v2}~\cite{reimers2019sentence}; details are given in Appendix~F. We specialize the reweighting map in \eqref{adaptweight} to the linear rule $\boldsymbol w_{\theta}=\theta \boldsymbol w$, where $\theta\in\mathbb{R}^{|\mathcal S|\times|\mathcal S|}$ satisfies Assumption~\ref{assumeaffine}. The resulting correction in \eqref{mixweight} is then obtained from the refined scenario-weight vector.}

\subsubsection{Preview variants}
We compare four preview variants.

\begin{itemize}
\item \textbf{Base.} Uses the fixed baseline preview $\hat{\pi}^{\mathrm{base}}_{0:47}$.

\item \textbf{Zero-shot \emph{LiFT}.} \x{Uses the initial scenario-weight vector $\boldsymbol w$ without online adaptation, yielding
$\hat{\pi}^{\mathrm{zero\text{-}shot}}_{0:47}
=
\hat{\pi}^{\mathrm{base}}_{0:47}
+
\sum_{s\in\mathcal S} w^s h^s_{0:47}$.}

\item \textbf{Online \emph{LiFT} (reset).} \x{Initializes $\theta=\boldsymbol I$ at the start of each day and updates $\boldsymbol w_\theta$ online via \eqref{adaptweight} and \eqref{eq:lift_update}. The resulting preview is
$\hat{\pi}_{t:\mathcal T(t)\mid t}
=
\hat{\pi}^{\mathrm{base}}_{t:\mathcal T(t)}
+
\sum_{s\in\mathcal S} w_\theta^s h^s_{t:\mathcal T(t)}$.}

\item \textbf{Online \emph{LiFT} (carry).} \x{Initializes $\theta=\boldsymbol I$ once and updates $\boldsymbol w_\theta$ online via \eqref{adaptweight} and \eqref{eq:lift_update} continuously across days, without resetting.}
\end{itemize}

On non-news days, all \emph{LiFT} corrections are skipped and the numerical baseline preview is used.

\subsection{Closed-loop control results}
\subsubsection{ESS arbitrage results}

Table~\ref{tab:econ_main} reports the economic reward 
$-\sum_{t=0}^{47}\pi_t u_t \Delta$,
evaluated under real prices over 2023--2024. 
All preview-based controllers use Chronos-T5 as the baseline numerical preview unless otherwise specified.

\vspace{-0.4cm}
\begin{table}[htbp]
\centering
\caption{Aggregated economic reward over 2023--2024 (AUD).}
\label{tab:econ_main}
\small
\begin{tabular}{lcc}
\toprule
\textbf{Controller} & \textbf{All days} & \textbf{Raw-news days} \\
\midrule
Off-OPT & $5.68\times10^7$ & $3.47\times10^7$ \\
Base-MPC & $5.61\times10^5$ & $3.64\times10^5$ \\
Zero-shot \emph{LiFT} & $1.13\times10^6$ & $9.38\times10^5$ \\
Online \emph{LiFT} (reset) & $2.59\times10^6$ & $2.39\times10^6$ \\
Online \emph{LiFT} (carry) & $2.59\times10^6$ & $2.40\times10^6$ \\
\bottomrule
\end{tabular}
\end{table}

The Off-OPT benchmark, with access to the full price trajectory, substantially outperforms all preview-based controllers, highlighting the value of accurate price information. Among implementable methods, \emph{LiFT} consistently improves over Base-MPC: Zero-shot \emph{LiFT} roughly doubles the reward, while online adaptation yields gains of over $4.6\times$ on all days and $6.5\times$ on raw-news days. The reset and carry variants perform nearly identically, indicating that most of the gain comes from within-day adaptation rather than \xiny{carrying the parameter across days}. Overall, these results demonstrate that \xiny{language-informed} preview refinement provides a substantial and robust economic benefit in closed-loop ESS operation.

\subsubsection{Operational behavior and constraint usage}
\x{Beyond economic reward, we also examine how often the closed-loop trajectory reaches the input and state limits. We find that \emph{LiFT} improves economic reward relative to Base-MPC without more frequent operation at the limits. In particular, the number of times the power limit is reached remains nearly unchanged across methods, while the number of times the SoC limit is reached decreases for Online \emph{LiFT} (reset). This suggests that the gains in Table~\ref{tab:econ_main} arise mainly from better timing of charge and discharge actions, rather than from more aggressive use of the operating limits. The reset and carry variants again exhibit nearly identical behavior.}

\subsubsection{Results with AR(48) baseline}

\begin{figure}[htbp]
    \centering
    \includegraphics[width=1.0\linewidth]{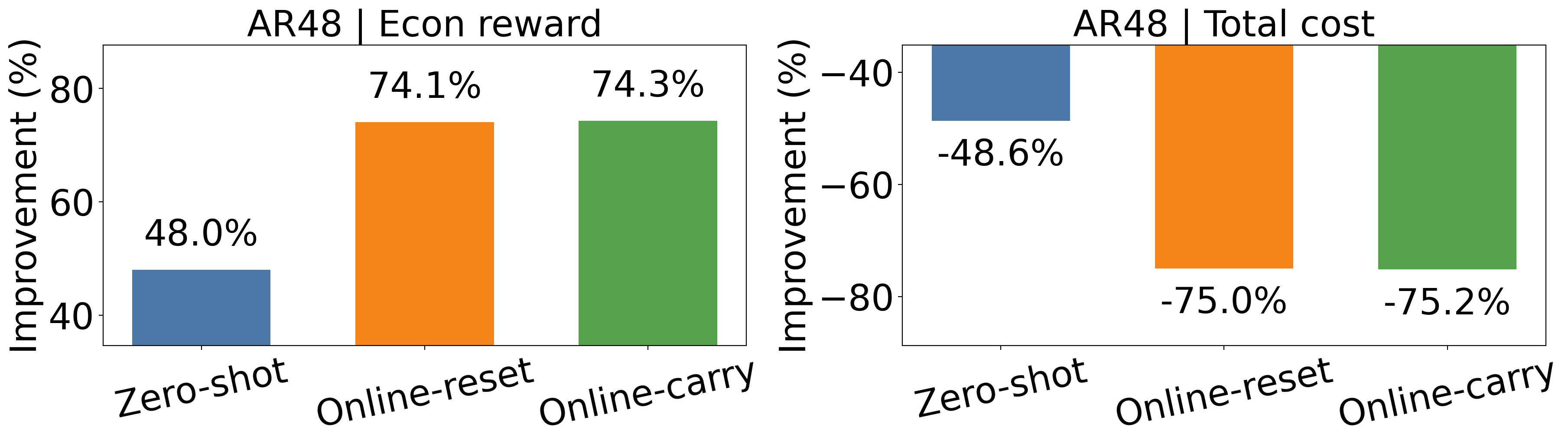}
    \vspace{-0.75cm}
    \caption{Percentage improvement over Base-MPC in cumulative economic reward and total cost for Zero-shot LiFT, Online-reset, and Online-carry, under AR(48) baseline over 2023--2024. Positive values indicate reward improvement, while negative values indicate cost reduction.}
    \label{fig:base_mpc_improvement_2x2}
\end{figure}

To evaluate robustness to weaker numerical previews, we repeat the experiments using AR(48) as the baseline predictor. Fig.~\ref{fig:base_mpc_improvement_2x2} reports the percentage improvement over Base-MPC in cumulative economic reward $-\sum_{t=0}^{47}\pi_t u_t \Delta$ and total cost $J^{\mathrm{MPC}}$.
Even with this weakest baseline, LiFT consistently delivers substantial gains. 
Zero-shot LiFT improves economic reward by $48.0\%$ and reduces total cost by $48.6\%$. 
With online adaptation, the gains further increase to over $74\%$ in reward and $75\%$ in cost reduction. Compared to the Chronos-based results, the relative improvements are smaller but remain significant, indicating that \emph{LiFT} provides robust benefits across different baseline qualities. 

\vspace{-0.45cm}
\vspace{-0.05cm}

\subsubsection{Effect of baseline quality on LiFT performance}
To further examine the role of the numerical baseline, we compare Zero-shot \emph{LiFT} built on top of Chronos (stronger) and AR(48) (weaker). For each day, forecasting accuracy is evaluated by comparing the predicted and realized 48-step intra-day price trajectories. The day-level mean absolute error (MAE) and root mean squared error (RMSE) are defined as
$\mathrm{MAE}=\frac{1}{48}\sum_{\tau=0}^{47}|\hat{\pi}_{\tau}-\pi_{\tau}|$
and
$\mathrm{RMSE}=\left(\frac{1}{48}\sum_{\tau=0}^{47}(\hat{\pi}_{\tau}-\pi_{\tau})^2\right)^{1/2}$.

\begin{figure}[htbp]
    \centering
    \includegraphics[width=1.0\linewidth]{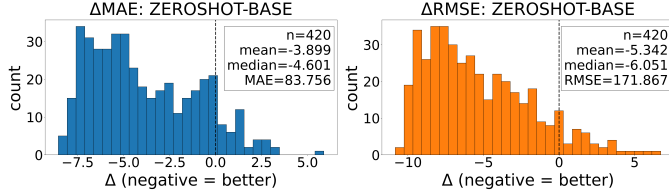}
    \vspace{-0.75cm}
    \caption{Histograms of day-level forecasting error differences between Zero-shot \emph{LiFT} and the Chronos baseline on 420 news days in 2023--2024.}
    \label{fig:delta_hist_mae_rmse_2x3}
\end{figure}

\begin{figure}[htbp]
    \centering
    \includegraphics[width=1.0\linewidth]{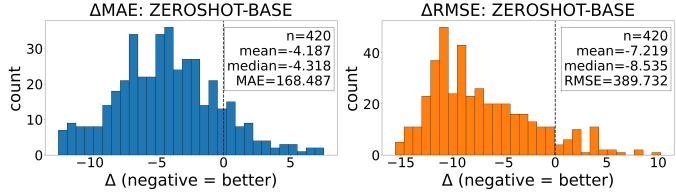}
    \vspace{-0.75cm}
    \caption{Histograms of day-level forecasting error differences between Zero-shot \emph{LiFT} and the AR(48) baseline on 420 news days in 2023--2024.}
    \label{fig:ardelta_hist_mae_rmse_2x3}
\end{figure}

Figs.~\ref{fig:delta_hist_mae_rmse_2x3} and~\ref{fig:ardelta_hist_mae_rmse_2x3} report the day-level error differences between Zero-shot \emph{LiFT} and the corresponding numerical baseline over the 420 news days in 2023--2024. Specifically, the horizontal axis shows
$\Delta\mathrm{MAE}:=\mathrm{MAE}_{\mathrm{Zero\text{-}shot}}-\mathrm{MAE}_{\mathrm{Base}}$
and
$\Delta\mathrm{RMSE}:=\mathrm{RMSE}_{\mathrm{Zero\text{-}shot}}-\mathrm{RMSE}_{\mathrm{Base}}$,
while the vertical axis gives the number of days in each bin. Negative values in the horizontal axis indicate that Zero-shot \emph{LiFT} improves on the baseline. Fig.~\ref{fig:delta_hist_mae_rmse_2x3} corresponds to Chronos, and Fig.~\ref{fig:ardelta_hist_mae_rmse_2x3} to AR(48). The mean and median reported in the figure are computed from these day-level differences, \x{whereas the displayed MAE and RMSE values are the averages of the corresponding day-level MAE and RMSE of Zero-shot \emph{LiFT} preview over the 420 news days.}

\x{The relative improvements are larger when \emph{LiFT} is built on the stronger baseline, Chronos. Relative to their respective baselines, Chronos-based \emph{LiFT} attains larger proportional reductions in MAE/RMSE than AR(48)-based \emph{LiFT}, suggesting that a stronger baseline enables more effective use of contextual information. Overall, the results show that \emph{LiFT} consistently improves performance across different numerical baselines, while the magnitude of the improvement depends on the quality of the baseline numerical preview.}

\vspace{-0.25cm}
\vspace{-0.35cm}

\vspace{-0.25cm}
\vspace{-0.35cm}

\section{Conclusion}
This paper proposes \emph{LiFT-MPC}, a language-in-the-loop feedback tuning framework for MPC with time-varying economic objectives. \emph{LiFT-MPC} augments a baseline numerical preview with a language-informed correction, admits delayed online updates via a control-performance loss, and provides performance and practical-stability guarantees. Numerical experiments on an energy-storage operation problem using real electricity prices and news context from New South Wales, Australia, show that \emph{LiFT} improves the closed-loop performance of BESS MPC. Even without online adaptation, \emph{LiFT} provides consistent preview-quality gains, while online adaptation further improves performance. Experiments with both Chronos and AR(48) baseline previews further demonstrate the benefits of the proposed \emph{LiFT-MPC} scheme.

\appendix

\noindent \textbf{A. Proof of Lemma~\ref{lem:quadratic_gap}:}
Let $V_T^{\mathrm{Off}}(x):=V_f(x)=x^\top P x$ and define
$Q_t^{\mathrm{Off}}(x,u):=\ell_t(x,u;\pi_t)+V_{t+1}^{\mathrm{Off}}(Ax+Bu),
V_t^{\mathrm{Off}}(x):=\min_u Q_t^{\mathrm{Off}}(x,u)$,
for $t=T-1,\ldots,0$. By dynamic programming, $V_0^{\mathrm{Off}}(x_0)=J^{\mathrm{Off\text{-}OPT}}$.
\xiny{From (\ref{eq:recursionansatz}) and (\ref{offoptlaw}), $Q_t^{\mathrm{Off}}(x,u)=u^\top H u+2b_t(x)^\top u+c_t(x)$ ,
where $H:=R+B^\top P B\succ0$, $b_t(x):=B^\top P A x+B^\top s_{t+1}^\star+\alpha_\pi \pi_t$, and $c_t(x)$ is independent of $u$. Completing the square gives $Q_t^{\mathrm{Off}}(x,u)
=
\bigl(u+H^{-1}b_t(x)\bigr)^\top H \bigl(u+H^{-1}b_t(x)\bigr)
+c_t(x)-b_t(x)^\top H^{-1}b_t(x)$.
By \eqref{offoptlaw}, $u_t^{\mathrm{Off}}(x)=-H^{-1}b_t(x)$, and hence $V_t^{\mathrm{Off}}(x)=c_t(x)-b_t(x)^\top H^{-1}b_t(x)$.
Therefore,}
{\small \begin{equation}\label{eq:dp_gap_identity}
Q_t^{\mathrm{Off}}(x,u)-V_t^{\mathrm{Off}}(x)
=
\bigl(u-u_t^{\mathrm{Off}}(x)\bigr)^\top H \bigl(u-u_t^{\mathrm{Off}}(x)\bigr).
\end{equation}}
At state $x_t^{\mathrm{MPC}}$, the minimizer of $Q_t^{\mathrm{Off}}(x_t,\cdot)$ is
\[
u_t^{\mathrm{Off}}(x_t^{\mathrm{MPC}})
=
-Kx_t^{\mathrm{MPC}}-H^{-1}\bigl(B^\top s_{t+1}^\star+\alpha_\pi\pi_t\bigr),
\]
whereas the implemented MPC input is
\[
u_t^{\mathrm{MPC}}(\theta_t)
=
-Kx_t^{\mathrm{MPC}}-H^{-1}\bigl(B^\top s_{t+1\mid t}(\theta_t)+\alpha_\pi\hat{\pi}_{t\mid t}(\theta_t)\bigr).
\]
Hence $u_t^{\mathrm{MPC}}(\theta_t)-u_t^{\mathrm{Off}}(x_t^{\mathrm{MPC}})
=
-H^{-1}\Bigl(
B^\top\bigl(s_{t+1\mid t}(\theta_t)-s_{t+1}^\star\bigr)
+
\alpha_\pi\bigl(\hat{\pi}_{t\mid t}(\theta_t)-\pi_t\bigr)
\Bigr)
=
\psi_t(\theta_t)$.
Applying \eqref{eq:dp_gap_identity} with
$x=x_t^{\mathrm{MPC}}$ and $u=u_t^{\mathrm{MPC}}$ gives
{\small \begin{equation}\label{recure2}
Q_t^{\mathrm{Off}}(x_t^{\mathrm{MPC}},u_t^{\mathrm{MPC}})
-
V_t^{\mathrm{Off}}(x_t^{\mathrm{MPC}})
=
\psi_t(\theta_t)^\top H\,\psi_t(\theta_t)=F_t.
\end{equation}}
Let $(x_t^{\mathrm{MPC}},u_t^{\mathrm{MPC}})$ denote the closed-loop trajectory induced by the preview-based MPC law (\ref{actmpc}). Evaluating the offline recursion along this trajectory and using
$x_{t+1}^{\mathrm{MPC}}=Ax_t^{\mathrm{MPC}}+Bu_t^{\mathrm{MPC}}$ gives
$Q_t^{\mathrm{Off}}(x_t^{\mathrm{MPC}},u_t^{\mathrm{MPC}})
=
\ell_t(x_t^{\mathrm{MPC}},u_t^{\mathrm{MPC}};\pi_t)
+
V_{t+1}^{\mathrm{Off}}(x_{t+1}^{\mathrm{MPC}})$.
Summing over $t=0,\ldots,T-1$ yields 
{\small \begin{align*}
 &\sum_{t=0}^{T-1}\Bigl(
Q_t^{\mathrm{Off}}(x_t^{\mathrm{MPC}},u_t^{\mathrm{MPC}})
-
V_t^{\mathrm{Off}}(x_t^{\mathrm{MPC}})
\Bigr)=\sum_{t=0}^{T-1}F_t\\
=&V_{T}^\mathrm{Off}-V_{0}^\mathrm{Off}+\sum_{t=0}^{T-1}\Bigl(
Q_t^{\mathrm{Off}}(x_t^{\mathrm{MPC}},u_t^{\mathrm{MPC}})
-
V_{t+1}^{\mathrm{Off}}(x_{t+1}^{\mathrm{MPC}})
\Bigr)\\
=&
J^{\mathrm{MPC}}(\theta_{0:T-1})-J^{\mathrm{Off\text{-}OPT}}.
\end{align*}}

\noindent \textbf{B. Proof of Lemma~\ref{lem:LD_bound_final}:}
Since $s_{i+1}^\star$ and $s^{\mathrm{real}}_{i+1\mid i}$ are independent of $\theta$, we have $\nabla_\theta \widehat\psi_i(\theta)=\nabla_\theta \psi_i(\theta)=R_i(\theta)$ and
$\widehat\psi_i(\theta)-\psi_i(\theta)
=
-H^{-1}B^\top\bigl(s^{\mathrm{real}}_{i+1\mid i}-s_{i+1}^\star\bigr)$.
Using $\nabla_\theta(\xi^\top H\xi)=2(\nabla_\theta\xi)^\top H\xi$ and the submultiplicative property,
{\small \begin{equation}\label{eq:LD_bound_mid_polished_appendix}
\begin{aligned}
&\bigl\|\nabla_\theta L_i(\theta)-\nabla_\theta F_i(\theta)\bigr\|\\
=&
2\bigl\|R_i(\theta)^\top H\bigl(\widehat\psi_i(\theta)-\psi_i(\theta)\bigr)\bigr\|\\
\le&
2\|R_i(\theta)\|\,\|H\|\,\|\widehat\psi_i(\theta)-\psi_i(\theta)\|\\
\le&
2\|R_i(\theta)\|\,\|H\|\,\|H^{-1}\|\,\|B\|\,
\bigl\|s^{\mathrm{real}}_{i+1\mid i}-s_{i+1}^\star\bigr\|\\
\le
&2\,\bar R\,\|H\|\,\|B\|\,\|H^{-1}\|\,
\bigl\|s^{\mathrm{real}}_{i+1\mid i}-s_{i+1}^\star\bigr\|.
\end{aligned}
\end{equation}}
Combining \eqref{eq:LD_bound_mid_polished_appendix} and \eqref{eq:s_trunc_bound_final_polished_appendix} yields \eqref{eq:LD_bound_explicit_polished}.

\noindent \textbf{C. Proof of Theorem~\ref{thm:regret_delay_ld}:}
By Lemma~\ref{lem:quadratic_gap}, $J^{\mathrm{MPC}}(\theta_{0:T-1})-J^{\mathrm{Off\text{-}OPT}}
=
\sum_{t=0}^{T-1}F_t(\theta_t),
J^{\mathrm{MPC}}(\theta)-J^{\mathrm{Off\text{-}OPT}}
=
\sum_{t=0}^{T-1}F_t(\theta)$,
for any fixed $\theta\in\Theta$. Hence $J^{\mathrm{MPC}}(\theta_{0:T-1})-J^{\mathrm{MPC}}(\theta^\star)
=
\sum_{t=0}^{T-1}\bigl(F_t(\theta_t)-F_t(\theta^\star)\bigr)$.
Under Assumption~\ref{assumeaffine}, $\hat{\pi}_{t\mid t}(\theta)$, $s_{t+1\mid t}(\theta)$, and thus $\psi_t(\theta)$ are affine in $\theta$. Therefore
$F_t(\theta)=\psi_t(\theta)^\top H\psi_t(\theta)$
is convex and differentiable, so $F_t(\theta_t)-F_t(\theta^\star)
\le
\nabla_\theta F_t(\theta_t)^\top(\theta_t-\theta^\star)$.
Summing over $t$ and adding/subtracting $\nabla_\theta L_t(\theta_t)$ gives
{\small \begin{equation}\label{eq:Ft_grad_bound_appendix}
\begin{aligned}
J^{\mathrm{MPC}}(\theta_{0:T-1})-J^{\mathrm{MPC}}(\theta^\star)
\le
&\sum_{t=0}^{T-1}\nabla_\theta L_t(\theta_t)^\top(\theta_t-\theta^\star)\\
&+\sum_{t=0}^{T-1}
\bigl(\nabla_\theta F_t(\theta_t)-\nabla_\theta L_t(\theta_t)\bigr)^\top\\
&(\theta_t-\theta^\star).
\end{aligned}
\end{equation}}
Using Cauchy--Schwarz, $\|\theta_t-\theta^\star\|\le D$, and Lemma~\ref{lem:LD_bound_final},
{\small \begin{equation}\label{eq:mismatch_bound_pf_appendix}
\sum_{t=0}^{T-1}
\bigl(\nabla_\theta F_t(\theta_t)-\nabla_\theta L_t(\theta_t)\bigr)^\top(\theta_t-\theta^\star)
\le
D\sum_{t=0}^{T-1}C_{\mathrm{LD}}\rho^{\mathcal T(t)-t}.
\end{equation}}

Now let $r_i:=\nabla_\theta L_i(\theta_i)$ and recall the delayed update
$\theta_{t+1}=\theta_t-\eta_t r_i$, with $r_i\equiv0$ for $t<k$. Then
$\sum_{t=0}^{T-1}\nabla_\theta L_t(\theta_t)^\top(\theta_t-\theta^\star)
=
\sum_{t=0}^{T-1}r_i^\top(\theta_t-\theta^\star)
+
\sum_{t=0}^{T-1}(r_t-r_i)^\top(\theta_t-\theta^\star)$.
\x{The first term is bounded by a standard OGD \ic{(online gradient descent)} telescoping argument; see, e.g., Eq.~(11) in \cite{wu2025instructmpc}. This gives
{\small \begin{equation}\label{eq:ogd_telescoping_pf_appendix}
\sum_{t=0}^{T-1}r_i^\top(\theta_t-\theta^\star)
\le
\frac{D^2}{2\eta_{T-1}}
+
\frac{G^2}{2}\sum_{t=0}^{T-1}\eta_t.
\end{equation}}}
\x{For the second term, a standard delayed-gradient estimate analogous to Eq.~(12) in \cite{wu2025instructmpc} yields
{\small \begin{equation}\label{eq:grad_diff_bound_loose_pf_appendix}
\sum_{t=0}^{T-1}(r_t-r_i)^\top(\theta_t-\theta^\star)
\le
(k-1)GD+(k-1)G^2\sum_{t=0}^{T-1}\eta_t.
\end{equation}}}

\noindent Combining \eqref{eq:Ft_grad_bound_appendix}--\eqref{eq:grad_diff_bound_loose_pf_appendix} yields \eqref{eq:regret_bound_general}.
For $\eta_t=\frac{D}{G\sqrt{t+1}}$,
$\frac{D^2}{2\eta_{T-1}}=\frac{DG\sqrt{T}}{2},
\frac{G^2}{2}\sum_{t=0}^{T-1}\eta_t
=
\frac{DG}{2}\sum_{t=0}^{T-1}\frac{1}{\sqrt{t+1}}
\le
DG\sqrt{T}$.
The first two terms in \eqref{eq:regret_bound_general} are bounded by $2DG\sqrt{T}$.

\noindent\textbf{D. Proof of Theorem~\ref{thm:practical_stab_short}:}
Let $e_t:=x_t-x_t^\star$.
Subtracting the state dynamics under~\eqref{eq:u_affine_split} and~\eqref{eq:u_star_affine_split} gives
$e_{t+1}
=
(A-BK)e_t+B(v_t-v_t^\star)
=
\Phi e_t+B\delta_t$.
Expanding the recursion yields
$e_t=\Phi^t e_0+\sum_{i=0}^{t-1}\Phi^{\,t-1-i}B\delta_i$.
Taking norms and using $\|\delta_i\|\le \bar\delta$ gives
$\|e_t\|
\le
c_\Phi\rho^t\|e_0\|
+
\sum_{i=0}^{t-1}c_\Phi\rho^{t-1-i}\|B\|\,\bar\delta$.
Therefore, $\|e_t\|
\le
c_\Phi\rho^t\|e_0\|
+
c_\Phi\|B\|\bar\delta\sum_{j=0}^{t-1}\rho^j
\le
c_\Phi\rho^t\|e_0\|
+
\frac{c_\Phi\|B\|}{1-\rho}\bar\delta$.

\noindent\textbf{E. \emph{LiFT} scenario library construction:}

\noindent\textbf{1) Residual clustering for historical trajectory prototypes $h^s_{0:47}$:}
For each day $d$ in 2018--2022 with available day-ahead news text, we compare the baseline numerical preview $\hat{\pi}^{\mathrm{base}}_{d,0:47}$ with the realized 48-step price trajectory $\pi_{d,0:47}$ and define the residual $\varepsilon_d:=\pi_{d,0:47}-\hat{\pi}^{\mathrm{base}}_{d,0:47}\in\mathbb R^{48}$.
We then remove its mean and normalize its magnitude:
$\tilde\varepsilon_d
=
\frac{\varepsilon_d-\bar\varepsilon_d\mathbf 1}{\|\varepsilon_d-\bar\varepsilon_d\mathbf 1\|_2},
\bar\varepsilon_d=\tfrac{1}{48}\mathbf 1^\top\varepsilon_d$.
\xiny{Each $\tilde\varepsilon_d\in\mathbb R^{48}$ is thus a normalized residual trajectory over the 48 half-hour intervals of one day. We apply K-means\cite{mcqueen1967some} clustering to the collection of normalized residual trajectories $\{\tilde\varepsilon_d\}$, so that trajectories with similar temporal shapes are grouped together. Each such group is referred to as a cluster. For each cluster $s\in\mathcal S$, its center is the empirical mean of the normalized residual trajectories assigned to that cluster, and therefore represents a typical residual-trajectory shape for that group. \xiny{We then rescale this cluster center by the median of the norms $\|\varepsilon_d\|_2$ over the trajectories assigned to that cluster, in order to obtain the historical trajectory prototype $h^s_{0:47}$.}}

\noindent\textbf{2) Scenario text grouping and semantic scenario text features $v^s$:}
For each cluster $s$, let $\mathcal C_s$ collect the associated day-ahead news texts. Each text $c\in\mathcal C_s$ is embedded by \texttt{sentence-transformers/all-MiniLM-L6-v2}~\cite{reimers2019sentence} as a normalized vector $\varphi(c)\in\mathbb R^{384}$. The normalized semantic feature vector is then defined as $v^s:=
\frac{\frac{1}{|\mathcal C_s|}\sum_{c\in\mathcal C_s}\varphi(c)}
{\left\|\frac{1}{|\mathcal C_s|}\sum_{c\in\mathcal C_s}\varphi(c)\right\|_2}$.

\noindent\textbf{\xinyy{F. The construction of the initial scenario-weighting vector $\boldsymbol w(c)\in\mathcal W^{|\mathcal S|}$:}} \xinyy{Given the context information $c$, we first compute its normalized semantic feature vector $\varphi(c)$ using the pretrained text encoder \texttt{sentence-transformers/all-MiniLM-L6-v2}~\cite{reimers2019sentence}. For each scenario $s\in\mathcal S$, let $v^s$ denote the normalized semantic feature vector. We then measure the relevance of scenario $s$ to the current context information through the cosine-similarity score}
\begin{equation}
\xinyy{r^s:=\cos\!\bigl(\varphi(c),v^s\bigr),\qquad s\in\mathcal S.}
\end{equation}
\xinyy{These similarity scores are converted into the initial scenario-weighting vector $\boldsymbol w(c)=\bigl[w^s(c)\bigr]_{s\in\mathcal S}$ by the temperature-softmax map}
\begin{equation}
\xinyy{w^s(c)=\frac{\exp(r^s/T)}{\sum_{j\in\mathcal S}\exp(r^j/T)},\qquad s\in\mathcal S,}
\end{equation}
\xinyy{where $T>0$ is a temperature parameter. By construction, $\boldsymbol w(c)\in\mathcal W^{|\mathcal S|}$, i.e., each entry is nonnegative and the entries sum to one. This vector provides the initial context-dependent scenario weighting, which is then adapted online through the reweighting map \eqref{adaptweight}.}

\bibliographystyle{IEEEtran}
\bibliography{CDCv0(03.31ddl)(Xinyi)} 

\end{document}